\documentclass[a4paper,11pt,english]{smfart}

\usepackage{amssymb}
\usepackage{url}
\usepackage[T1]{fontenc}
\RequirePackage{calrsfs}
\DeclareSymbolFont{rsfscript}{OMS}{rsfs}{m}{b}
\DeclareSymbolFontAlphabet{\mathrsfs}{rsfscript}
\usepackage[utopia,expert]{mathdesign}
\usepackage{lscape}
\usepackage{array, boldline, makecell, booktabs}

\usepackage{enumitem}

\usepackage{todonotes}

\usepackage[leftbars]{changebar}

\usepackage{palatino}
\usepackage{rotating}
\usepackage{graphicx}

\usepackage{tikz-cd}

\usepackage{enumerate}

\usepackage{amscd}

\usepackage{color}
\definecolor{shadecolor}{gray}{0.90}

\input xypic
\xyoption{all}
\xyoption{arc}

\makeindex

\newtheorem{theo}{Theorem}[section]

\newtheorem{lem}[theo]{Lemma}

\def\equat{\refstepcounter{theo}\begin{equation}}
\def\endequat{\end{equation}}

    \def\CM{{\mathbb{C}}}

    \def\QM{{\mathbb{Q}}}

\def\Ab{{\mathbf A}}    \def\AC{{\mathcal{A}}}
    \def\BC{{\mathcal{B}}}
    \def\CC{{\mathcal{C}}}

\def\Gb{{\mathbf G}}    
    
    \def\IC{{\mathcal{I}}}
    \def\JC{{\mathcal{J}}}
    
\def\Lb{{\mathbf L}}    \def\LC{{\mathcal{L}}}

\def\Pb{{\mathbf P}}

    \def\SC{{\mathcal{S}}}

    \def\ZC{{\mathcal{Z}}}

\def\Zrm{{\mathrm{Z}}}

\def\a{\alpha}
\def\b{\beta}

\def\l{\lambda}

\def\Sig{\Sigma}
\def\th{\theta}

\def\mub{{\boldsymbol{\mu}}}

\DeclareMathOperator{\Gal}{{\mathrm{Gal}}}

\DeclareMathOperator{\Id}{{\mathrm{Id}}}

\DeclareMathOperator{\Proj}{{\mathrm{Proj}}}

\def\to{\rightarrow}

\def\DS{\displaystyle}

\def\lexp#1#2{\kern\scriptspace\vphantom{#2}^{#1}\kern-\scriptspace#2}
\def\le{\hspace{0.1em}\mathop{\leqslant}\nolimits\hspace{0.1em}}
\def\ge{\hspace{0.1em}\mathop{\geqslant}\nolimits\hspace{0.1em}}

\mathchardef\inferieur="321E
\mathchardef\superieur="321F

\def\eqna{\begin{eqnarray*}}
\def\endeqna{\end{eqnarray*}}

\def\itemth#1{\item[${\mathrm{(#1)}}$]}

\catcode`\@=11
\long\def\@car#1#2\@nil{#1}
\long\def\@first#1#2{#1}
\long\def\@second#1#2{#2}
\long\def\ifempty#1{\expandafter\ifx\@car#1@\@nil @\@empty
  \expandafter\@first\else\expandafter\@second\fi}
\catcode`\@=12

\def\GL{{\mathrm{GL}}}

\theoremstyle{remark}

\theoremstyle{plain}

\def\BIL{LR}
\def\GAUCHE{L}
\def\CAR{CAR}
\def\FAM{FAM}

\def\xyinj{\ar@{^{(}->}}
\def\xysur{\ar@{->>}}

\bigskip

\def\trespetitespace{\vphantom{$\DS{\frac{\DS{A}}{\DS{A}}}$}}

\makeatletter
\def\hlinewd#1{%
\noalign{\ifnum0=`}\fi\hrule \@height #1 %
\futurelet\reserved@a\@xhline}
\makeatother

\newlength\epaisLigne

\usepackage{multirow,multicol}

\makeatletter
\def\hlinewd#1{%
\noalign{\ifnum0=`}\fi\hrule \@height #1 %
\futurelet\reserved@a\@xhline}
\makeatother

\usepackage{multirow}

\def\trespetitespace{\vphantom{$\DS{\frac{\DS{A}}{\DS{A}}}$}}

\usepackage{lscape}
\def\zero{{\boldsymbol{0}}}

\def\GL{\operatorname{\Gb\Lb}\nolimits}

\begin{document}

\title{A dodecic surface with 320 cusps}

\author{{\sc C\'edric Bonnaf\'e}}
\address{IMAG, Universit\'e de Montpellier, CNRS, Montpellier, France}

\makeatletter
\email{cedric.bonnafe@umontpellier.fr}
\makeatother

\date{\today}

\thanks{The author is partly supported by the ANR 
(Project No ANR-24-CE40-3389, GRAW)}

\begin{abstract}
We construct a degree $12$ homogeneous invariant of the complex reflection group
$G_{29}$ (in Shephard-Todd's notation) whose associated
surface has 320 singularities of type $A_2$, thereby improving previous records for dodecic
surfaces.
\end{abstract}

\maketitle
\pagestyle{myheadings}
\markboth{\sc C. Bonnaf\'e}{A dodecic surface with 320 cusps}

\bigskip

For a type $X$ of isolated surface singularity, and for $d \ge 1$, let $\mu_X(d)$
denote the maximal number of singularities of type $X$ a surface of degree $d$
in $\Pb^3(\CM)$ might have.
Determining $\mu_X(d)$ is a classical problem. The exact value is known only in a few
cases (for instance for $X=A_1$ and $d \le 6$) but there have been many works
for trying to give upper and lower bounds for $\mu_X(d)$. For quotient singularities,
and for $d$ even or $d \ge 14$, the best-known upper bounds are given by Miyaoka~\cite{miyaoka}.

Lower bounds are obtained by constructing explicit examples of surfaces of degree $d$
with many singularities.
For type $A$ singularities and for any $d$, lower bounds have been obtained by several
authors~\cite[\!\dots]{chmutov, labs these, labs Aj, escudero1,escudero2,escudero3}. In small degrees and in type
$A_1$, some exceptional
examples give better lower bounds (see the works of Togliatti~\cite{togliatti}, Barth~\cite{barth},
Endra\ss~\cite{endrass}, Labs~\cite{labs these, labs}, Sarti~\cite{sarti 0,sarti 1,sarti},...).

For small degrees and singularities of type $A_2$ (which will be called {\it cusps} in this paper),
Labs~\cite[Theo.~7.1]{labs these} proved that $\mu_{A_2}(6) \ge 35$ while
Borisov and Gunnels~\cite{borisov} proved that $\mu_{A_2}(8) \ge 84$:
these are the best-known lower bounds.
Our aim in this paper is to improve the best-known lower bound for $\mu_{A_2}(12)$,
obtained by Escudero~\cite{escudero2}, who
constructed a dodecic surface with $301$ cusps:

\bigskip

\noindent{\bf Theorem.} {\it There exists an homogeneous invariant of degree $12$ of the complex reflection
group $G_{29}$ (in Shephard-Todd's notation~\cite{shephard todd}) whose zero
set is a surface with exactly $320$ cusps, which form a single $G_{29}$-orbit.}

\bigskip

Together with Miyaoka's upper bound, our result says that
$$320 \le \mu_{A_2}(12) \le 363.$$
It must be said that the results of Escudero are much more general, since he
got lower bounds for all degrees divisible by $3$ and $A_j$-singularities
for any $j \ge 2$. For instance, he
proved that
$$96k^2(4k-1)+14k-1 \le \mu_{A_2}(12k) \le 3k(12k-1)^2$$
(here, the upper bound is again due to Miyaoka).
Using the classical trick consisting in lifting a surface of degree $d$ in $\Pb^3(\CM)$ to
a surface of degree $kd$ through the morphism $\Pb^3(\CM) \to \Pb^3(\CM)$,
$[x:y:z:t] \mapsto [x^k:y^k:z^k:t^k]$
allows to construct, thanks to our surface, a surface of degree $12k$
with $320k^3$ cusps: this gives the lower bound $\mu_{A_2}(12k) \ge 320k^3$, which is better than
Escudero's one only for $k=1$.

We also investigate several other singular dodecic surfaces defined by a fundamental
invariant of $G_{29}$, some of them might be of interest for algebraic geometers.
In this list, we retrieve a singular dodecic surface with $160$ singularities
of type $D_4$ which was already constructed by the author~\cite{bonnafe}, and still
gives the bestknown lower bound for $\mu_{D_4}(12)$.

\bigskip
%
%

\noindent{\bf Notation.} We set $V=\CM^4$ and we denote by $(x,y,z,t)$ the dual basis of
the canonical basis of $\CM^4$: the algebra $\CM[V]$ of polynomial functions on $V$
is the polynomial algebra $\CM[x,y,z,t]$. We identify $\GL_\CM(V)$ with $\GL_4(\CM)$
and the projective space $\Pb(V)$ with
$\Pb^3(\CM)$.

If $m$ is a monomial in $x$, $y$, $z$ and $t$, we denote by $\Sig_4(m)$
the sum of all the monomials obtained from $m$ by permutation of these four variables.
For instance,
$$\Sig_4(xy)=xy+xz+xt+yz+yt+zt\qquad\text{and}\qquad \Sig_4(xyzt)=xyzt.$$

\section{The complex reflection group ${\boldsymbol{G_{29}}}$}

\medskip

\subsection{Definition}
Let $i \in \CM$ denote a square root of $-1$ and let
$$s_1=
\begin{pmatrix}
. & 1 & . & . \\
1 & . & . & . \\
. & . & 1 & . \\
. & . & . & 1 \\
\end{pmatrix},\qquad s_2=
\begin{pmatrix}
1 & . & . & . \\
. & . & 1 & . \\
. & 1 & . & . \\
. & . & . & 1 \\
\end{pmatrix},$$
$$s_3=\begin{pmatrix}
. & -i & . & . \\
i & . & . & . \\
. & . & 1 & . \\
. & . & . & 1
\end{pmatrix}\quad\text{and}\quad
s_4=\frac{1}{2}
\begin{pmatrix}
\hphantom{-}1 & -1 & -1 & -1 \\
-1 & \hphantom{-}1 & -1 & -1 \\
-1 & -1 & \hphantom{-}1 & -1 \\
-1 & -1 & -1 & \hphantom{-}1
\end{pmatrix}.
$$
The matrices $s_1$, $s_2$, $s_3$ and $s_4$ are reflections of order $2$ and they generate a finite
subgroup of $\GL_4(\CM)$, which can be taken as a model for the complex reflection group
denoted by $G_{29}$ in Shephard-Todd's classification~\cite{shephard todd}. So we set
$$G_{29}=\langle s_1,s_2,s_3,s_4\rangle \subset \GL_4(\QM[i]) \subset \GL_4(\CM).$$
Note that $G_{29}$ is irreducible and that
\equat\label{eq:order}
|G_{29}|=7\,680\qquad\text{and}\qquad \Zrm(G_{29}) =\mub_4 \Id_V.
\endequat

\bigskip

\subsection{Invariants}
Note also that $G_{29}$ is stable under the complex conjugacy (because $\overline{s}_3=s_1s_3s_1$).
Since $\Gal(\QM[i]/\QM)$ is generated by the complex conjugacy, a theorem of
Marin-Michel~\cite{marin-michel} implies that one can choose a family of fundamental
invariants in $\QM[x,y,z,t]$. Moreover, the subgroup $G_{29} \cap \GL_4(\QM)$ is the rational
reflection group denoted by $G(2,1,4)$ in Shephard-Todd's classification (it is isomorphic
to the Weyl group of type $B_4$). In particular, any $G_{29}$-invariant polynomial
is a linear combination of elements of the form $\Sig_4(m)$, where $m$ is a monomial
in $x^2$, $y^2$, $z^2$ and $t^2$. We set
$$
\begin{cases}
f_1=\Sig_4(x^4)-6\,\Sig_4(x^2y^2),\\
f_2=\Sig_4(x^8) + 4\,\Sig_4(x^6y^2) + 6\,\Sig_4(x^4y^4)
- 20\,\Sig_4(x^4y^2z^2) + 152\,x^2y^2z^2t^2,\\
f_3=\Sig_4(x^8y^2z^2) - \Sig_4(x^6y^4z^2) + 2\,\Sig_4(x^6y^2z^2t^2) \\
\hskip6cm - 2\,\Sig_4(x^4y^4z^4) + 2\,\Sig_4(x^4y^4z^2t^2).
\end{cases}$$
Then $f_1$, $f_2$, $f_3$ are homogeneous $G_{29}$-invariant polynomials of respective degrees $4$, $8$ and $12$
and there exists an homogeneous invariant $f_4$ of degree $20$ such that
$$\CM[V]^{G_{29}}=\CM[f_1,f_2,f_3,f_4].$$
Note that $f_3$ is, up to a scalar, the unique fundamental invariant of degree $12$ whose degree in $x$ is $\le 8$.
Also, $f_2$ is, up to a scalar, the Hessian of $f_1$.

\bigskip

\begin{lem}\label{lem:irreductible}
Any fundamental invariant of degree $12$ of $G_{29}$ is irreducible.
\end{lem}

\bigskip

\begin{proof}
Let $F$ be a fundamental invariant of degree $12$ of $G_{29}$ and let
$f$ be an irreducible divisor of $F$. Note that $f$ is necessarily homogeneous.
We set
$$G=\{g \in G_{29}~|~g(f) \in \CM^\times f\}.$$
Write $f_\sharp=\prod_{g \in [G_{29}/G]} g(f)$, where $[G_{29}/G]$
is a set of representatives of the cosets in $G_{29}/G$. Then $g(f_\sharp) \in \CM^\times f_\sharp$
for all $g \in G_{29}$. Note that $f_\sharp$ divides $F$ because
$F$ is $G_{29}$-invariant. Note also that $f$ is $G'$-invariant (where $G'$
denotes the derived subgroup of $G$) because the map $G \to \CM^\times$,
$g \mapsto g(f)/f$ is a linear character of $G$.

Now, let $\th : G_{29} \to \CM^\times$, $g \mapsto g(f_\sharp)/f_\sharp$: it is a linear character
of $G_{29}$.
There are only two linear characters of $G_{29}$, namely the trivial one
and the restriction of the determinant. If $\th$ is the restriction of the determinant,
then it follows for instance from~\cite[Theo.~9.19]{lehrer taylor} that $f_\sharp$
has degree bigger than the number of reflecting hyperplanes of $G_{29}$, which is equal to $40$.
This is impossible because $f_\sharp$ divides $F$.

This shows that $f_\sharp$ is $G_{29}$-invariant. Three cases may occur:
\begin{itemize}
\item[$\bullet$] If $\deg(f_\sharp)=4$, then $f_1$ divides $F$, so $F/f_1$ is homogeneous of degree $8$ and $G_{29}$-invariant,
so it is of the form
$\a f_2+\b f_1^2$ for some $\a$, $\b \in \CM$. This is impossible since  this would
give an algebraic relation between $F$, $f_1$ and $f_2$.

\item[$\bullet$] If $\deg(f_\sharp)=8$, then $F/f_\sharp$ is an homogeneous $G_{29}$-invariant divisor
of degree $4$ of $F$. We conclude that it is impossible as in the previous case.

\item[$\bullet$] If $\deg(f_\sharp)=12$, then this shows that
$$F= \kappa \prod_{g \in [G_{29}/G]} g(f)$$
for some $\kappa \in \CM^\times$. Write $d=\deg(f)$ and $r=|G_{29}/G|$. Then
$dr=12$ and $f$ is $G'$-invariant. But one can easily check with {\sc Magma}~\cite{magma} that, for any subgroup $G$
of $G_{29}$ of index $r \in \{2,3,4,6,12\}$, the derived subgroup of $G$
has no non-zero homogeneous invariant of degree $12/r$. This shows that
$d=12$ and $r=1$, i.e. that $F=\kappa f$, as expected.
\end{itemize}
The proof of the lemma is complete.
\end{proof}

\section{The family of invariant dodecics}

\medskip

If $f \in \CM[V]$ is homogeneous, we denote by $\ZC(f)$ the scheme
$\Proj\bigl(\CM[V]/\langle f \rangle\bigr)$.
Any reduced, irreducible surface defined by an homogeneous invariant of
degree $12$ of $G_{29}$ is of the form $\SC_{12}^{\l,\mu}=\ZC(f_3+\l f_2f_1 + \mu f_1^3)$
for some $(\l,\mu) \in \Ab^2(\CM)$. Lemma~\ref{lem:irreductible} says that the converse holds,
that is,
\equat\label{eq:s12}
\text{\it $\SC_{12}^{\l,\mu}$ is reduced and irreducible}
\endequat
for any $(\l,\mu) \in \Ab^2(\CM)$.

\bigskip

\subsection{Singular invariant dodecics}
The subset $\CC$ of $\Ab^2(\CM)$ formed by the elements $(\l,\mu)$ such that $\SC_{12}^{\l,\mu}$
is singular is a closed subset of $\Ab^2(\CM)$. The computation of $\CC$
was too long for our computer: adapting slightly the algorithm of~\cite[\S{3}]{bonnafe},
we computed with {\sc Magma} the subset $\CC_0$ of $\CC$
defined as the set of elements $(\l,\mu) \in \Ab^2(\CM)$ such that
$\SC_{12}^{\l,\mu}$ has a singular point having its last
coordinate equal to $0$. Then $\CC_0$ is of pure dimension $1$ and is the
union of $8$ irreducible curves:
\begin{itemize}
\item[$\bullet$] Six lines $\LC_1$, $\LC_2$, $\LC_3$, $\LC_4$, $\LC_5^+$ and $\LC_5^-$ defined by the equations
$$
\begin{cases}
(\LC_1) &\hskip0.7cm \mu=0,\\
(\LC_2) &\hskip0.7cm \l=0,\\
(\LC_3) &\hskip0.7cm \l+\mu=0,\\
(\LC_4) &\hskip0.7cm \l - 15\,\mu - \frac{1}{45}=0,\\
(\LC_5^+) &\hskip0.7cm \l -(4+2i)\,\mu - \frac{3+i}{320} =0,\\
(\LC_5^-) &\hskip0.7cm \l -(4-2i)\,\mu - \frac{3-i}{320} =0,\hphantom{AAAAAAAA}\\
\end{cases}
$$

\item[$\bullet$] A cubic curve $\AC$ defined by the equation
$$20\,480\,\l^3 - 256\,\l^2 + \l + \mu=0.$$

\item[$\bullet$] A sextic curve $\BC$ defined by
$$1\,342\,177\,280\,\l^6 - 100\,663\,296\,\l^5 + 3\,014\,656\,\l^4 - 3\,538\,944\,\l^3\mu\hskip2cm$$
$$\hskip2cm - 45\,056\,\l^3 + 73\,728\,\l^2\mu + 336\,\l^2 - 288\,\l\mu - \l - 432\,\mu^2 - \mu=0.$$
\end{itemize}
We do not know whether $\CC_0$ and $\CC$ are equal.

\medskip

\subsection{320 cusps}
The first seven irreducible components are isomorphic to $\Ab^1(\CM)$ and $\BC$
has two singular points $(\l^\pm,\mu^\pm)$, where
$$
\l^\pm=\frac{3\pm\sqrt{3}}{384}\qquad\text{and}\qquad
\mu^\pm= \frac{-5\pm 3\sqrt{3}}{6\,912}.
$$
We set
$$\SC_{12}^\pm=\ZC(f_3+\l^\pm f_2f_1 + \mu^\pm f_1^3).$$
The main result of this paper is the following: the proof is obtained through a
straightforward {\sc Magma}
computation, which takes only a few minutes on a standard laptop.

\bigskip

\begin{theo}\label{theo:main}
The surfaces $\SC_{12}^+$ and $\SC_{12}^-$ are reduced, irreducible, of degree $12$
and admit $320$ cusps and no other singular point. These cusps form a single $G_{29}$-orbit.
\end{theo}

\bigskip

Unfortunately, even though the surfaces $\SC_{12}^\pm$ are defined over $\QM(\sqrt{3})$,
none of their singular point is real, so drawing it with {\sc Surfer}~\cite{surfer}
does not lead to a beautiful picture. We do not know whether this surface can be defined
over $\QM$.

\bigskip

\subsection{Some other singular dodecics}
The following can also be checked by a computer calculation with {\sc Magma}:
\begin{itemize}
\itemth{1} For generic $(\l,\mu)$ in $\LC_1$, the surface $\SC_{12}^{\l,\mu}$ has $20$ singularities of type $T_{4,4,4}$\footnote{A singularity is said to be of type $T_{4,4,4}$ if it is equivalent to the singularity at the point $0$
of the affine surface $\{(a,b,c) \in \Ab^3(\CM)~|~abc+a^4+b^4+c^4=0\}$. Such a singularity has multiplicity $3$, Milnor number
$11$ and Tjurina number $10$.}.

\itemth{2} For generic $(\l,\mu)$ in $\LC_2$, the surface $\SC_{12}^{\l,\mu}$ has $120$ singularities of type $A_1$.

\itemth{3} For generic $(\l,\mu)$ in $\LC_3$, the surface $\SC_{12}^{\l,\mu}$ has $40$ singularities of type $A_1$.

\itemth{4} For generic $(\l,\mu)$ in $\LC_4$, the surface $\SC_{12}^{\l,\mu}$ has $160$ singularities of type $A_1$.

\itemth{5} For generic $(\l,\mu)$ in $\LC_5^\pm$, the surface $\SC_{12}^{\l,\mu}$ has $80$ singularities of type $A_1$.

\itemth{6} For generic $(\l,\mu)$ in $\AC$, the surface $\SC_{12}^{\l,\mu}$ has $480$ singularities of type $A_1$.

\itemth{7} For generic $(\l,\mu)$ in $\BC$, the surface $\SC_{12}^{\l,\mu}$ has $320$ singularities of type $A_1$.
\end{itemize}
Theorem~\ref{theo:main} shows for instance that, for $(\l,\mu)$ running in $\BC$, the $320$ singularities
of type $A_1$ degenerate to $A_2$-singularities when $(\l,\mu)$ reaches a singular point of $\BC$.
Note also that the two singular points $(\l^\pm,\mu^\pm)$ of $\BC$ do not belong to any
other irreducible component of $\CC_0$.

Since $\BC$ is the unique singular irreducible component of $\CC_0$, the singular locus
of $\CC_0$ consists of $(\l^\pm,\mu^\pm)$ and the points lying on at least two
irreducible components of $\CC_0$. There are 39 such intersection points.
We have computed with {\sc Magma} the singularities of the 39 associated
dodecic surfaces and we have checked the following facts:
\begin{itemize}
\itemth{a} If $(\l_0,\mu_0)$ belongs to only two irreducible components $\IC$ and $\JC$ of $\CC_0$
and if $\IC$ and $\JC$ intersect transversely at $(\l_0,\mu_0)$, then the dodecic surface $\SC_{12}^{\l_0,\mu_0}$
cumulates the singularities ``coming from $\IC$'' and those ``coming from $\JC$''.
There are 33 singular points of $\CC_0$ satisfying this property.
For instance:
\begin{itemize}
\item[$\bullet$] The point $p_{1,4}=(1/45,0)$ belongs only to $\LC_1$ and $\LC_4$, so $\SC_{12}^{p_{1,4}}$
has 20 singularities of type $T_{4,4,4}$ and 160 singularities of type $A_1$. These are its
only singular points: they form two $G_{29}$-orbits.



\smallskip

\item[$\bullet$] The point $p_5^\pm=(-(-1\pm i)/320,-(3\pm i)/1600)$ belongs only to $\LC_5^\pm$ and $\AC$, and $\LC_5^\pm$ and $\AC$
intersect transversely at this point, so $\SC_{12}^{p_5^\pm}$ has 560 singularities
of type $A_1$, and no other singular point: they form two $G_{29}$-orbits, of respective cardinality
$80$ and $480$.
\end{itemize}

\smallskip

\itemth{b} The point $\clubsuit=(1/40,1/5400)$ belongs only to $\LC_4$ and $\BC$, and $\LC_4$ and $\BC$ intersect
at $\clubsuit$ with multiplicity $2$. The surface $\SC_{12}^{\clubsuit}$ has 160 singularities of type $A_3$
and no other singular point: they form a single $G_{29}$-orbit.

\smallskip

\itemth{c} The point $\diamondsuit=(1/240,-13/10\,800)$ belongs only to $\LC_4$ and $\AC$, and $\LC_4$ and $\AC$ intersect
at $\diamondsuit$ with multiplicity $3$. The surface $\SC_{12}^{\diamondsuit}$ has 160 singularities of type $D_4$
and no other singular point: they form a single $G_{29}$-orbit. This is the best-known lower bound for
$\mu_{D_4}(12)$: this surface was already discovered by the author~\cite[Ex.~5.5(3)]{bonnafe}.

\smallskip

\itemth{d} The point $\heartsuit=(1/64,0)$ belongs only to $\LC_4$ and $\BC$, and $\LC_4$ and $\BC$ intersect
at $\heartsuit$ with multiplicity $4$. The surface $\SC_{12}^{\heartsuit}$ has 20 singular points,
which form a single $G_{29}$-orbit:
they have multiplicity $4$, Milnor number $27$ and Tjurina number $26$ and their projective tangent cone is smooth.
Using the software {\sc Singular}, we checked that this singularity is not in Arnold's list: let us denote
this type of singularity by $X_\heartsuit$.

\smallskip

\itemth{e} The point $\spadesuit^\pm=((3\pm i)/640,(-7\pm i)/6\,400)$ belongs to $\LC_5^\pm$, $\AC$ and $\BC$
and to no other irreducible component of $\CC_0$. The surface $\SC_{12}^{\spadesuit^\pm}$ has
$80$
singularities of type $T_{4,4,4}$ and no other singular points. They form a single $G_{29}$-orbit.

\smallskip

\def\zero{{\boldsymbol{0}}}
\itemth{f} The point $\zero=(0,0)$ belongs to $\LC_1$, $\LC_2$, $\LC_3$, $\AC$ and $\BC$ and does not belong
to $\LC_4$ and $\LC_5^\pm$. The singular locus of the surface $\SC_{12}^\zero$ is the union
of 30 lines (in particular, it is of dimension $1$): these lines form a single $G_{29}$ orbit,
namely the orbit of the line defined by $z=t=0$. This surface was also mentioned in~\cite[Ex.~5.5(4)]{bonnafe}.

\end{itemize}
\medskip

\begin{centerline}
{\begin{tabular}{@{{{\vrule width 2pt}\,\,}}c@{{\,\,\,{\vrule width 1pt}}\,\,}c@{{\,\,\,{\vrule width 2pt}}}}
\hlinewd{2pt}
\trespetitespace Surface & Singularities \\
\hlinewd{1pt}
\trespetitespace $\SC_{12}^\pm$ & 320 $A_2$ \\
\hline
\trespetitespace $\SC_{12}^{\clubsuit}$ & 160 $A_3$ \\
\hline
\trespetitespace $\SC_{12}^{\diamondsuit}$ & 160 $D_4$ \\
\hline
\trespetitespace $\SC_{12}^{\heartsuit}$ & 20 $X_\heartsuit$ \\
\hline
\trespetitespace $\SC_{12}^{\spadesuit^\pm}$ & 80 $T_{4,4,4}$ \\
\hlinewd{2pt}
\end{tabular}}
\end{centerline}
%
%
\medskip

\centerline{$\SC_{12}^{p_{1,4}}$\includegraphics[width=5cm]{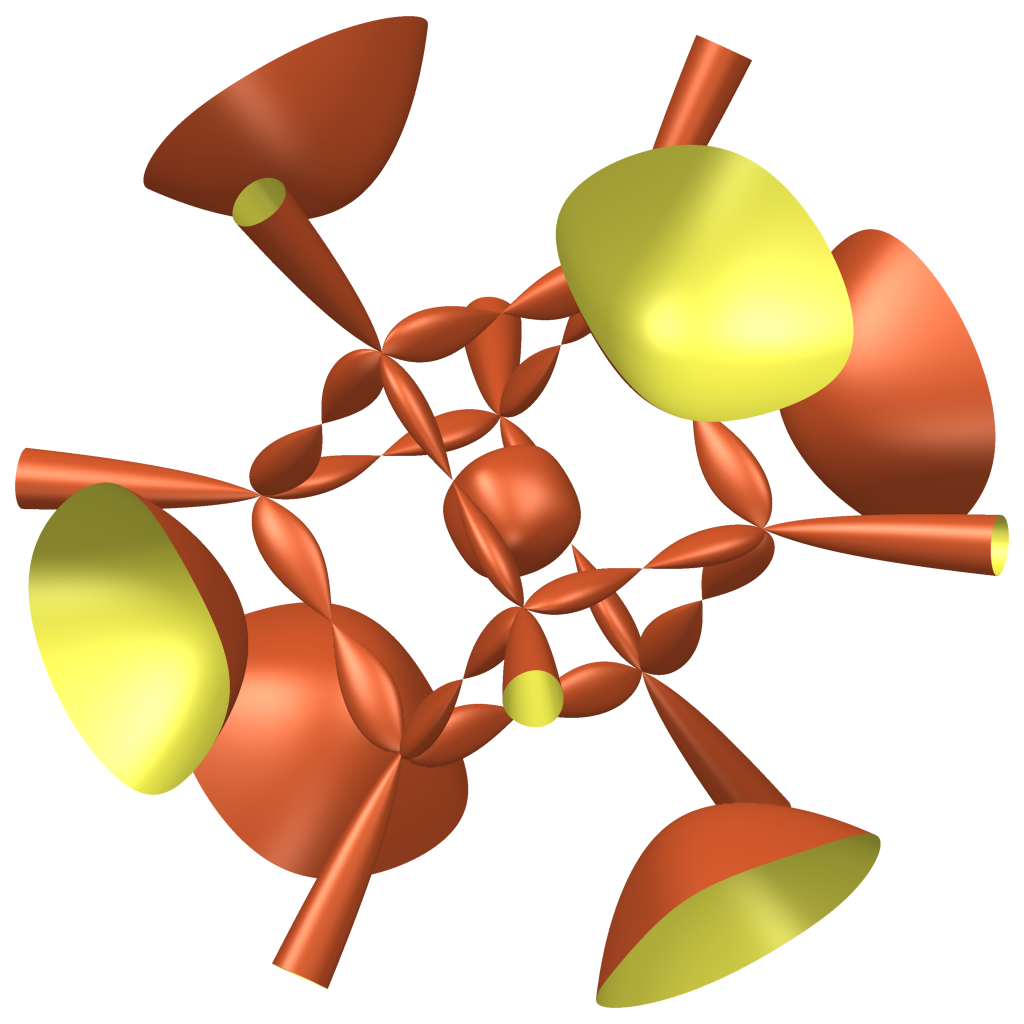} \hskip1.5cm
\includegraphics[width=5cm]{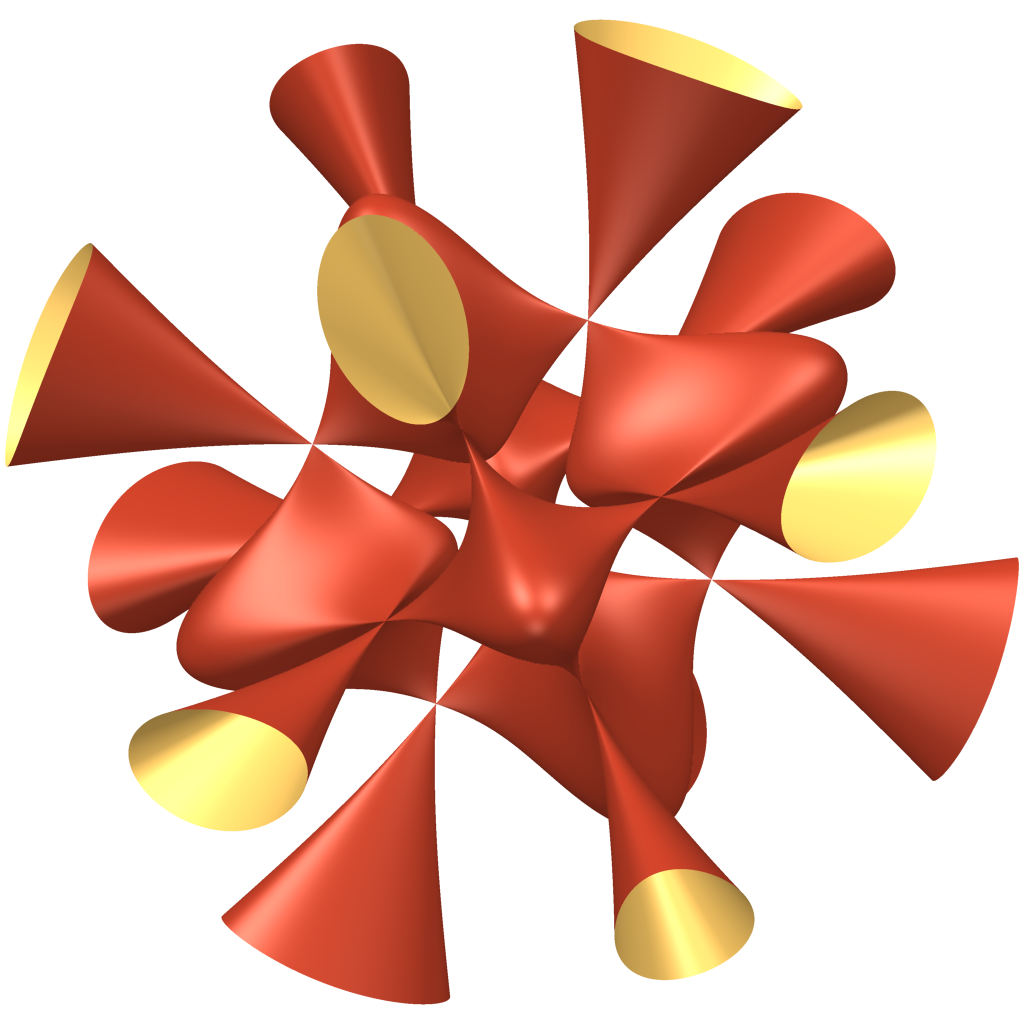}$\SC_{12}^{\diamondsuit}$}

%


\centerline{\includegraphics[width=5cm]{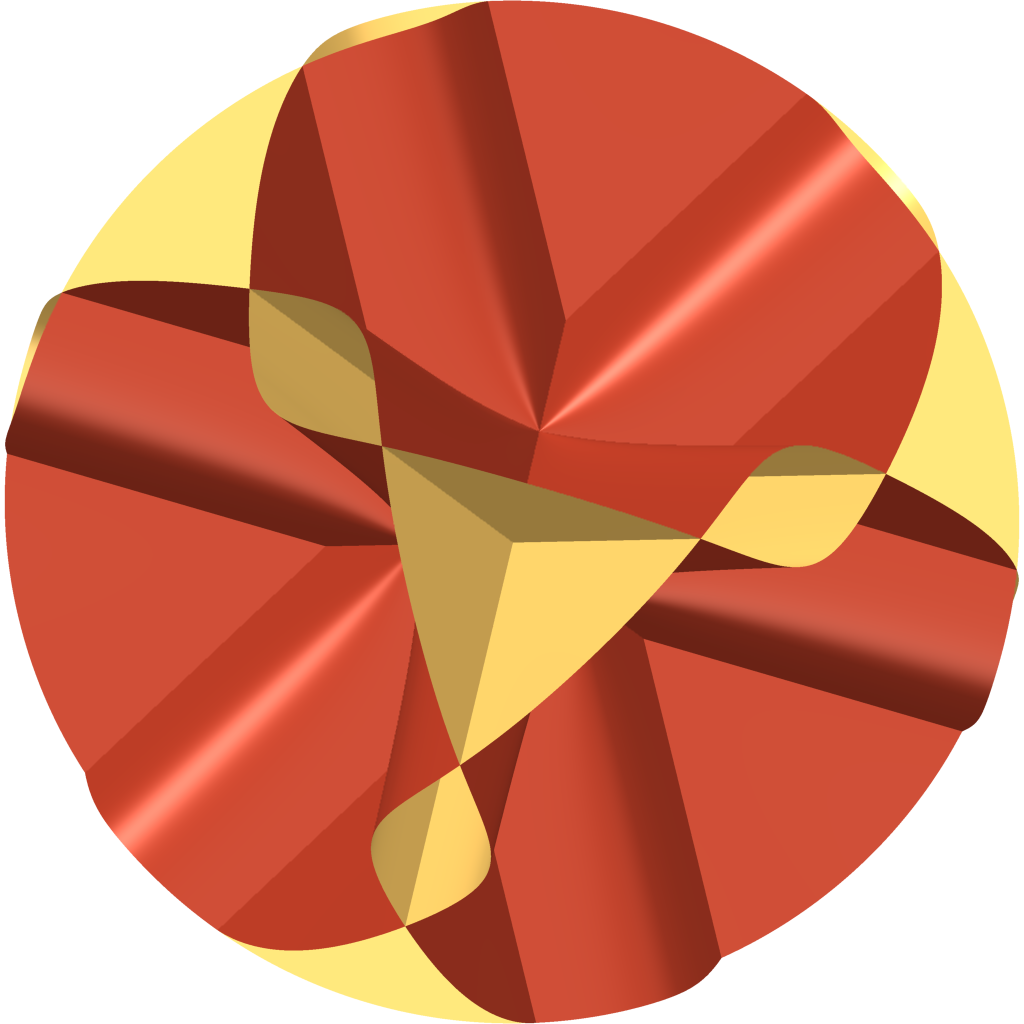}}

\smallskip

\centerline{$\SC_{12}^\zero$}

\end{document}